\documentclass[11pt]{amsart}
\usepackage{epsfig}
\usepackage{subfigure}
\usepackage{psfrag}

\usepackage{graphicx}
\usepackage{amscd}
\usepackage{amsmath}
\usepackage{amsfonts}
\usepackage{yfonts}
\usepackage{psfrag}
\usepackage{amssymb}
\usepackage{verbatim}
\usepackage{comment}

\newtheorem{theorem}{Theorem}[section]
\theoremstyle{plain}

\newtheorem{corollary}[theorem]{Corollary}

\newtheorem{defi}[theorem]{Definition}

\newtheorem{lemma}[theorem]{Lemma}

\newtheorem{prop}[theorem]{Proposition}
\newtheorem{remark}[theorem]{Remark}

\def\inrad{{\rm inr}}
\def\Comp{{\rm Comp}}

\def\diam{{\rm diam\,}}
\def\dist{{\rm dist}}

\def\vphi{\varphi}

\def\Ok{{\mathcal O}}
\def\Mk{{\mathcal M}}

\def\Th{\theta}

\def\wtil{\widetilde}

\def\im{{\rm Im}}
\def\re{{\rm Re}}

\newcommand{\lam}{\lambda}

\newcommand{\gam}{\gamma}
\newcommand{\om}{\omega}

\newcommand{\Sig}{\Sigma}

\newcommand{\sig}{\sigma}

\newcommand{\R}{{\mathbb R}}

\newcommand{\disc}{{\mathbb D}}
\newcommand{\Int}{{\rm int}}

\def\N{{\mathbb N}}
\newcommand{\C}{{\mathbb C}}
\def\D{{\mathbb D}}

\newcommand{\Nat}{{\mathbb N}}

\def\Bk{{\mathcal B}}
\def\Dk{{\mathcal D}}

\def\Tk{{\mathcal T}}
\def\ba{{\bf a}}
\def\bB{{\bf b}}
\def\bw{{\bf w}}
\def\half{\frac{1}{2}}

\newcommand{\Fk}{{\mathcal F}}

\newcommand{\es}{\emptyset}

\def\ov{\overline}

\begin{document}

\thispagestyle{empty}

\title[Quasisymmetric conjugacy]{Quasisymmetric conjugacy between quadratic dynamics and iterated
function systems}

\author{Kemal Ilgar Ero\u{g}lu}
\address{Kemal Ilgar Ero\u{g}lu, Department of Mathematics and Statistics P.O. Box 35 (MaD) FI-40014 University of Jyv\"askyl\"a Finland}
\email{kieroglu@hotmail.com}

\author{Steffen Rohde}
\address{Steffen Rohde, Box 354350, Department of Mathematics,
University of Washington, Seattle WA 98195}
\email{rohde@math.washington.edu}

\author{Boris Solomyak}
\address{Boris Solomyak, Box 354350, Department of Mathematics,
University of Washington, Seattle WA 98195}
\email{solomyak@math.washington.edu}

\thanks{ 
\indent Research of Ero\u{g}lu and Solomyak was
partially supported by NSF grants \#DMS-0355187 and \#DMS-0654408. Research of Rohde was
supported in part by NSF Grants DMS-0501726
and DMS-0800968.}

\begin{abstract}
We consider linear iterated function systems (IFS) with a constant contraction ratio
in the plane for which the ``overlap set'' $\Ok$ is finite, and which are
``invertible'' on the attractor $A$, in the sense that there is a continuous surjection $q: A\to A$ whose inverse branches
are the contractions of the IFS. The overlap set is the critical set in the sense that $q$ is not a local homeomorphism precisely
at $\Ok$. We suppose also that there is a rational function $p$ with the Julia set $J$ such that $(A,q)$ and $(J,p)$ are
conjugate. We prove that if $A$ has bounded turning and $p$ has no parabolic cycles, then the conjugacy is quasisymmetric.
This result is applied to some specific examples including an uncountable family. Our main focus is on the family of IFS $\{\lambda z,\lambda z+1\}$ where 
$\lambda$ is a complex parameter in the unit disk, such that its attractor $A_\lam$ is a dendrite, which happens whenever 
$\Ok$ is a singleton. C.~Bandt observed that a simple
modification of such an IFS (without changing the attractor) is invertible and gives rise to a quadratic-like map $q_\lam$ on $A_\lam$.
If the IFS is post-critically finite, then a result of A. Kameyama shows that there is a quadratic map $p_c(z)=z^2+c$,
with the Julia set $J_c$ such that $(A_\lam,q_\lam)$ and $(J_c,p_c)$ are conjugate. We prove that this conjugacy is
quasisymmetric and obtain partial results in the general (not post-critically finite) case. 
\end{abstract}

\date{\today}

\maketitle

\thispagestyle{empty}

\section{Introduction}

For $\lam$ in the open unit disk $\D$
consider the compact set $A_\lam\subset \C$ given by
\begin{equation} \label{eq-al}
A_\lam = \Bigl\{\sum_{n=0}^\infty a_n \lam^n:\ a_n \in
\{0,1\}\Bigr\}\,.
\end{equation}
It is the attractor for the iterated function system (IFS) $\{F_0,F_1\}$,
where $F_0(z)=\lam z,\ F_1(z)=\lam z + 1$,
that is, $A_\lam$ is the unique non-empty
compact set such that $A_\lam = F_0(A_\lam) \cup F_1(A_\lam)=
\lam A_\lam \cup (1+\lam A_\lam)$, see \cite{hutch}.
The connectedness locus
$$
\Mk:= \{\lam\in \D:\ A_\lam \mbox{\ is connected}\}
$$
was introduced in
\cite{BH} and studied by several authors, see \cite{bou2,bandt,solxu,solo}
in particular.
There is the same dichotomy as for Julia sets $J_c$ of
quadratic maps $p_c(z)=
z^2+c$: the attractor $A_\lam$ is either connected or totally disconnected.
It is well-known that $A_\lam$ is connected if and only if
$
\lam A_\lam \cap (\lam A_\lam + 1) \ne \es.
$
It follows that $\Mk$ has a characterization as the set of zeros of
$\{-1,0,1\}$ power series:
$$
\Mk = \Bigl\{\lam\in \disc:\ \exists\,b_k \in \{-1,0,1\},\ k\ge 1,\ \mbox{such that}\ 1 + \sum_{k=1}^\infty b_k\lam^k =0\Bigr\}.
$$

Figure 1 shows the part of $\Mk$ in
$\{z:\ \re(z)\ge 0,\,\im(z)> 0,\,|z| < 1/\sqrt{2}\}$ (the set $\Mk$ is
symmetric with respect to both axes, and all $\lam\in \disc$, with
$|\lam|\ge 1/\sqrt{2}$, are known to be in $\Mk$).
Bousch \cite{bou2} proved that $\Mk$ is
connected and locally connected.

\begin{figure}[lt]
\epsfig{figure=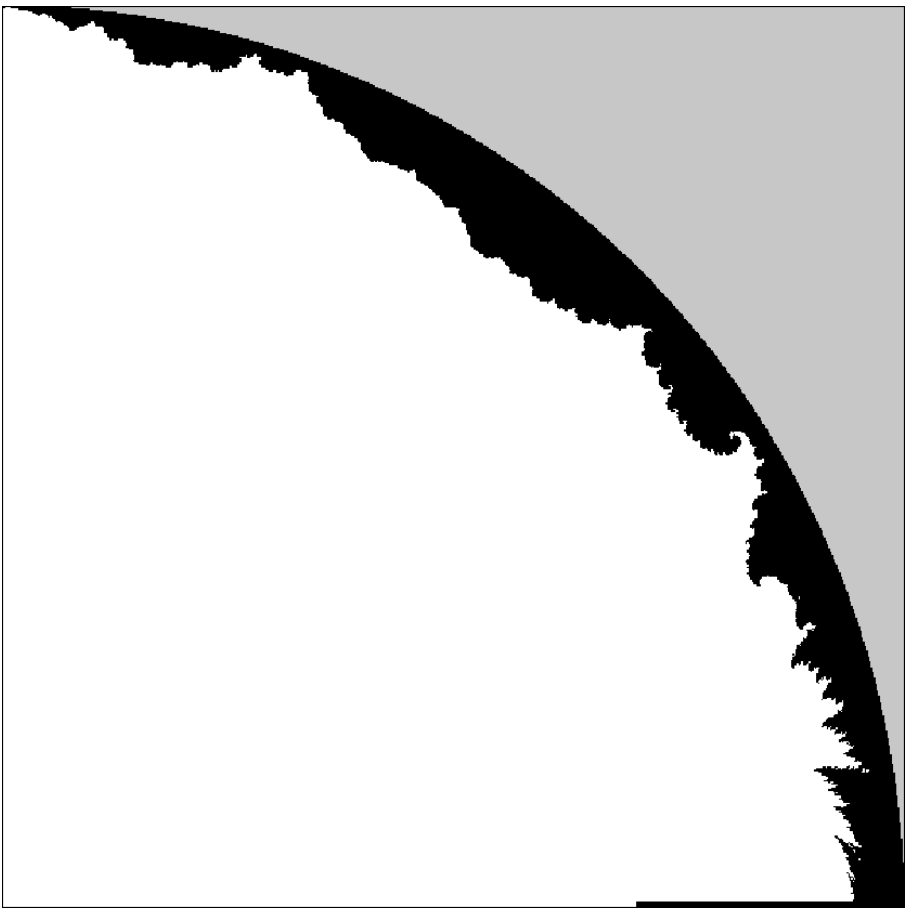,height=12cm}
\caption{The connectedness locus $\Mk$ (figure made by C. Bandt)}
\end{figure}

Here we consider the following subset of $\Mk$:
$$
\Tk:= \{\lam\in \D:\ \lam A_\lam \cap (\lam A_\lam + 1) = \{o_\lam\} \mbox{\ is a singleton}\}.
$$
It is known \cite{bake} (and not hard to see) that $\lam\in \Tk$ if
and only if $A_\lam$ is
a {\em dendrite}, that is,
a connected, locally connected, nowhere dense compact set in
the plane with connected complement. 
(Note that a connected self-similar set
is necessarily locally connected, see \cite{hata}.)

The set $A_\lam$ is invariant under the involution $s(z) = - z+(1-\lam)^{-1}$.
(On the symbolic level this map just switches the coefficients 0 and 1 in
the power series representation.) 
It follows that $\lam A_\lam \cap (\lam A_\lam + 1)$ is invariant under $s$ , hence for $\lam\in \Tk$ we have
$o_\lam = \frac{1}{2(1-\lam)}$, the single common point of $\lam A_\lam$ and $\lam A_\lam + 1$.
Bandt \cite{bandt} observed that it is possible
to define quadratic-like dynamics on $A_\lam$ for $\lam\in \Tk$.
Indeed, $A_\lam$ is also the
attractor of the IFS $\{F_0, F_1 \circ s\}$, and these two contractions are
the inverse branches of a function
\begin{equation} \label{eq-ql}
q_\lam(z)=
        \left\{ \begin{array}{ll} F_0^{-1}(z) = \frac{z}{\lam}, & \mbox{if}
                 \ z \in F_0(A_\lam)=\lam A_\lam;\\
          s\circ F_1^{-1}(z) =        
        \frac{1-z}{\lam} + \frac{1}{(1-\lam)}, &
                   \mbox{if} \ z \in F_1(A_\lam)=\lam A_\lam+1;
\end{array} \right.
\end{equation}
Note that $q_\lam(o_\lam)$ is well-defined, hence $q_\lam$ is continuous, it maps
$A_\lam$ onto itself 2-to-1, except at the single critical point $o_\lam$;
at all other points it is a local homeomorphism.
This is an example of what Kameyama \cite{kam1} called an {\em invertible IFS}.
We say that the system $(A_\lam,q_\lam)$ is {\em 
post-critically finite} (p.c.f.) if $o_\lam$ is
strictly preperiodic under $q_\lam$.
Bandt \cite{bandt} observed that by results of Kameyama \cite{kam1,kam2},
if $(A_\lam,q_\lam)$ is p.c.f., then it 
is topologically conjugate with 
$(J_c,p_c)$ for a certain $c=c(\lam)$.
The sets $A$ and $J$ in Figure 2 and 3 are ``more similar'' than mere topological
equivalence would explain. It appears visually  that they are quasiconformal
images of each other. In fact, the attempt to justify this visual observation 
was the starting point of our work and led to the following.

\begin{theorem} \label{th-quad1}
Suppose $\lam\in \Tk$ is such that $(A_\lam,q_\lam)$ is post-critically
finite. Then the conjugacy with $(J_c,p_c)$ is quasisymmetric and
extends to a quasiconformal homeomorphism of the Riemann sphere.
\end{theorem}

We do not know if the same holds for all $\lam\in \Tk$, without the
p.c.f.\ assumption, but we have been able to prove it for an uncountable set $\Tk_0\subset \Tk$, see (\ref{def-T0}) and
Corollary~\ref{cor-new}. 

The first ingredient of the proof is establishing the
{\em bounded turning property}, defined in Section \ref{BTPI},
closely related to the concept of a
{\em John domain}.
We prove the bounded turning property for a large class of
self-similar fractals, which includes p.c.f.\ $A_\lam$ for $\lam\in \Tk$.
This part does not depend on the invertibility
of the IFS.

Our proof of quasisymmetry is similar in
spirit to the proof of McMullen and Sullivan \cite{MS}
of quasiconformality of conjugacies between hyperbolic rational maps, 
and its extension
to non-hyperbolic settings as in \cite{PR}. See \cite{HP} for an 
exposition of some of these ideas in a general setting.
Our proof applies not only to the sets $A_\lam$, but to other
examples as well, such as the Sierpi\'nski gasket and the map
$p(z) = z^2 - \frac{16}{27z}$, see Section 4.
In the proof we use  a recent result of Bandt and Rao \cite{BR} that
a planar self-similar IFS with a connected attractor and a
finite ``overlap set'' is non-recurrent.

\medskip

{\bf Examples.} The simplest example is $\lam = 1/2$ and $c=-2$, which is
a well-known pair of quasisymmetrically
conjugate maps, since $q_\lam$ is just the tent-map for $\lam = 1/2$.
In \cite{solo} there are many examples of $\lam$ for which
it is rigorously proven that $A_\lam$ is a dendrite (including an
uncountable set which we denote $\Tk_0$ in this paper),
and many more can be found numerically, using pictures and
a version of  Bandt's algorithm from \cite{bandt}.
In the figures below we show two examples to which our theorem applies.
For $\lam\in \Tk$ we indicate the (unique) $\{-1,0,1\}$ power series  vanishing at $\lam$,
by the sequence of its coefficients, e.g.,
$+(-+++--)$; here $+$ and $-$ correspond to $\pm 1$'s and
parentheses indicate a period. For the Julia set $J_c$ we also indicate the external parameter angle. See Section 5 for an explanation
of these examples.

\begin{figure}[ht]
\begin{center}
$\begin{array}{cc}
\multicolumn{1}{l}{} &
        \multicolumn{1}{l}{} \\
\epsfxsize=2.in
\epsffile{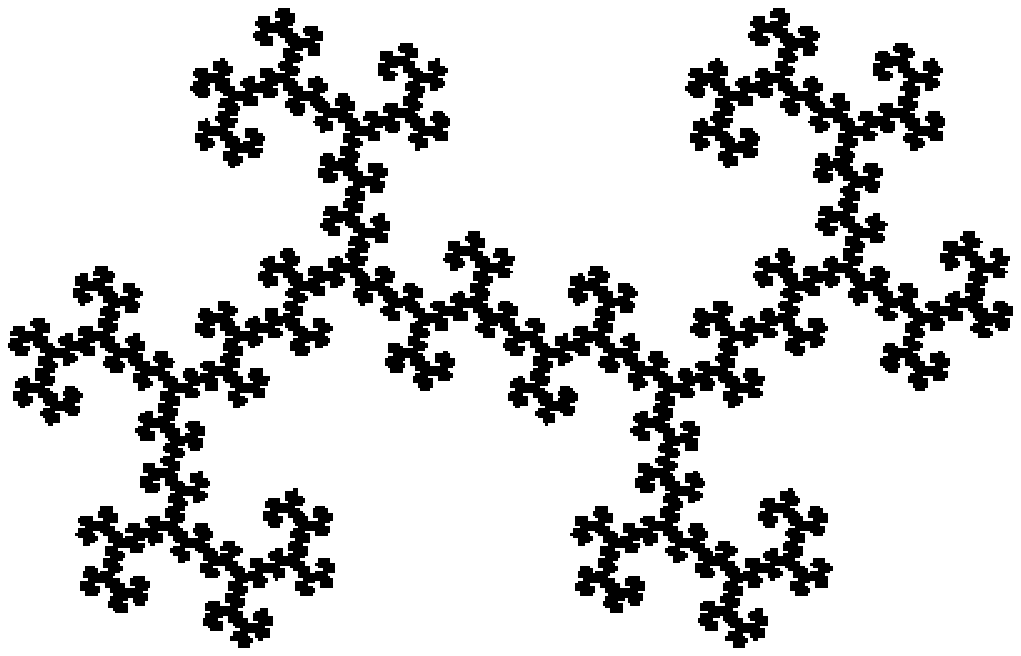} \ \ & \ \
        \epsfxsize=2.5in
      \epsffile{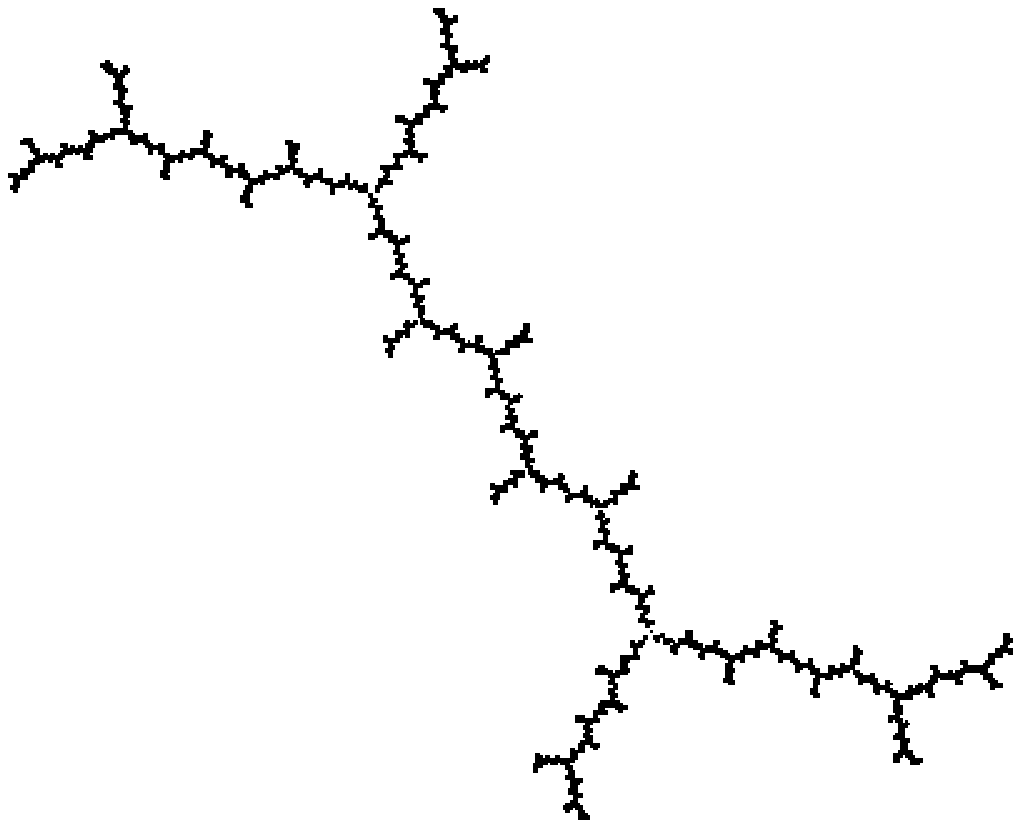} \\[1.cm]
 \mbox{Self-similar set $A_\lam$ } \ \ & \ \ \mbox{Julia set $J_c$}
\end{array}$
\end{center}
\caption{$\lam\approx .3668760 + .5202594i$, the power series $+(-+++--)$,
the kneading sequence $1\ov{100}$,
$c\approx -.155788+1.11222i$, external angle $\frac{3}{14}$.}
\end{figure}

\begin{figure}[ht]
\begin{center}
$\begin{array}{cc}
\multicolumn{1}{l}{} &
        \multicolumn{1}{l}{} \\
\epsfxsize=2.3in
\epsffile{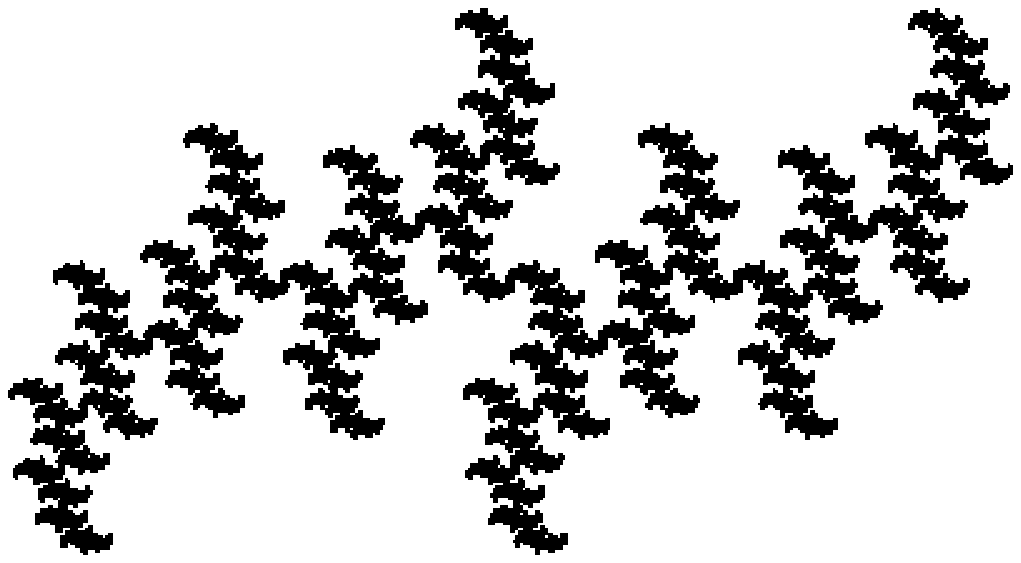} \ \ & \ \
        \epsfxsize=2.6in
        \epsffile{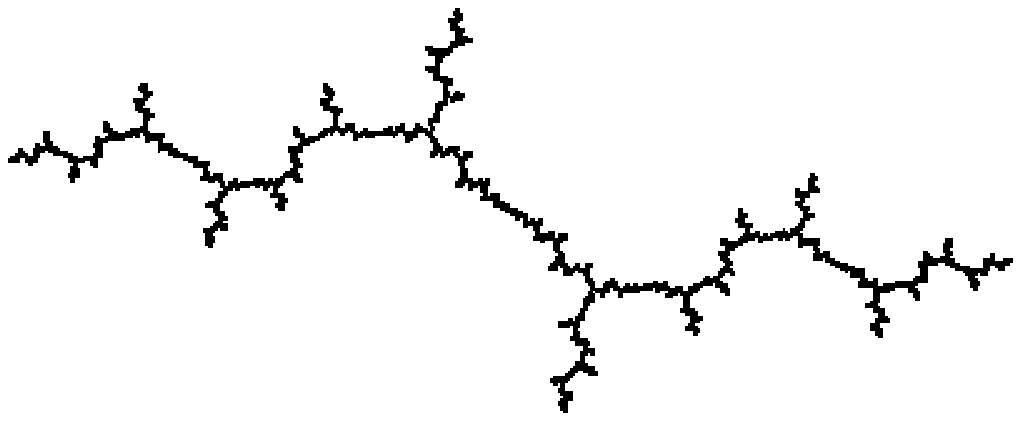} \\[1.cm]
\mbox{Self-similar set $A_\lam$} \ \ & \ \ \mbox{Julia set $J_c$}
\end{array}$
\end{center}
\caption{$\lam\approx .595744 + .254426i$, the power series
$+--(+)$, the kneading sequence $101\ov{0}$,
$c\approx -1.29636+.441852i$, external angle $\frac{3}{8}$.}
\end{figure}

\medskip

The rest of the paper is organized as follows.
In the next section we prove the bounded turning property for the sets $A_\lam$, using a ``dynamical'' argument.
In Section 3 we prove, in a more general setting, that if there is a conjugacy between piecewise linear and rational dynamics,
then it is quasisymmetric. Sections 2 and 3 yield Theorem~\ref{th-quad1}. Section 4 is devoted to examples of rational (non-polynomial)
maps for which our theorem applies. Sections 5 and 6 deal with various extensions, in particular,
to the non-recurrent case. 

\medskip

{\bf Notation.} Throughout this paper, $F$ and $f$ denotes affine maps (similitudes),
$p$ denotes a polynomial or rational map, $p_c(z)=z^2+c$, $q$ is a piecewise linear
expanding map, and $B(x,r)$ is the closed disk of radius $r$ centered at $x.$

\section{Bounded Turning Property I}\label{BTPI}

A connected, locally connected set $X\subset \C$ has {\em bounded turning} if there exists
$L$ such that for every $z_0,z_1 \in X$ there is a connected set $[z_0,z_1]
\subset X$ satisfying
\begin{equation} \label{eq-bt}
\diam[z_0,z_1] \le L |z_0-z_1|.
\end{equation}

\begin{theorem} \label{th-bt1} Let $\lam \in \Tk$, $A_\lam$ be given by (\ref{eq-al}) and $q_\lam$ by (\ref{eq-ql}), and suppose
that $(A_\lam,q_\lam)$ is post-critically finite. Then $A_\lam$ has the bounded turning property.
\end{theorem}

\noindent
{\em Proof.}
Let $q=q_\lam$,
$A=A_\lam$, $A_j = \lam A + j$ for $j=0,1$, so that $A = A_0 \cup A_1$. Recall that $o_\lam$ is the only ``critical point,''
the intersection of $A_0$ and $A_1$. The assumption that  $(A,q)$ is post-critically finite means that for some $\ell \ge 1$ and
$p\ge 1$, we have $q^{\ell+p}(o_\lam) = q^\ell(o_\lam)$, and $q^n(o_\lam)  \ne o_\lam$ for $n\ge 1$. 

Consider $z_0\ne z_1$ in $A$, with the goal of proving (\ref{eq-bt}). We can assume that $z_0, z_1$ do not lie in the same
subset 
$A_j$, $j=0,1$. Indeed, otherwise we can pass to $z_j' = q(z_j)$, $j=0,1$, using the fact  that $q$ is an expanding similitude on
$A_j$, so the ``local geometry'' is preserved. Eventually, the points will be separated by $o_\lam$. So we may suppose that $z_j \in A_j$,
$j=0,1$, and none of the points is $o_\lam$.

Let $\delta>0$, to be determined later. We have
\begin{equation} \label{eq-eta}
\eta(\delta):= \max_{j=0,1}\dist(A_j, A_{1-j}\setminus \Int B(o_\lam,\delta))>0,
\end{equation}
in view of $A_0 \cap A_1 = \{o_\lam\}$. Therefore,
if $r=\max_{j=0,1} |z_j-o_\lam| \ge \delta$, then 
(\ref{eq-bt}) holds with $[z_0, z_1]= A$ and $L=\diam A/\eta(\delta),$
and we are done. We will reduce the general case to this case by using 
the self-similarity of $A$ at $q^{\ell}(o_\lam)$.
Let $k\ge 0$ be the smallest integer such that
\begin{equation} \label{eq-tuchka}
r|\lam|^{-pk}=\max_{j=0,1} |z_j-o_\lam|\, |\lam|^{-pk} \ge \delta.
\end{equation}
The idea is to ``enlarge'' the picture by applying $q^{-\ell} q^{\ell+p k}$ for a suitable
branch of the inverse. In order to avoid dealing with different branches, we use the 
symmetry of $A$ about $o_\lam$ and consider $z'_1=s(z_1).$
We choose $\delta$ so small that $o_\lam\notin q^{j}(B(o_\lam,\delta))$ for all $1\le j\le \ell+p$.
Then $q^{\ell+p k}|_{A_0\cap B(o_\lam,r)}$ is an expanding similitude. There are two branches
of $q^{-\ell}$ that map $q^{\ell}(o_\lam)$ back to $o_\lam.$ Denote $Q$ the one that has
$Q(A)\subset A_0$. Then $R=Q  q^{\ell+p k}$ is an expanding similitude on $A_0\cap B(o_\lam,r)$ which
fixes $o_\lam$ and maps $z_0,z'_1$ to points $\zeta_0, \zeta'_1.$
Set $\zeta_1 = s(\zeta'_1),$ then 
$$\max_{j=0,1} |\zeta_j - o_\lam| \ge \delta$$
in view of (\ref{eq-tuchka}).
Hence $$|\zeta_0-\zeta_1| \ge \eta(\delta)$$ by (\ref{eq-eta}) and therefore
$$|z_0-z_1| =|R^{-1}\zeta_0-sR^{-1}\zeta_1'| = |R^{-1}\zeta_0 - R^{-1}\zeta_1| =|\lam|^{pk}|\zeta_0-\zeta_1|\geq \eta(\delta) |\lam|^{p k}.$$
Here we used that $R^{-1}$ extends as an affine linear map to the plane which fixes $o_\lambda$, so it commutes with the involution $s$.
Since $R^{-1}(A)$ contains $z_0, z'_1$
and $o_\lam$,
the set $[z_0,z_1] =  R^{-1}(A)\cup s(R^{-1}(A))$ has
$$\diam[z_0,z_1] \le  2|\lam|^{pk}\diam A \leq
 \frac{2\diam A}{\eta(\delta)} |z_0-z_1|,
$$
as desired .\qed

\section{Quasisymmetry of the conjugacy}

Here we show that the homeomorphism $\vphi$ conjugating
the quadratic polynomial $p_c$ on its Julia set $J_c$ to the piecewise linear
map $q_\lam$ on $A_\lam$ is quasisymmetric on $J_c$, in the setting of Theorem~\ref{th-quad1}.
In fact, the proof works in much greater generality, so we start with
a precise description of our setting.

We suppose that $A$ is the attractor of an IFS $\{f_j\}_{j=1}^m$ 
in the complex plane, with $m\ge 2$, where $f_j$ are similarities, all having
the same contraction ratio $|\lam|<1$ (so we can write $f_j(z) = |\lam|
e^{i\Th_j}z + d_j$ for some $\Th_j\in [0,2\pi)$ and $d_j\in \C$, if 
$f_j$ is orientation-preserving, and $f_j(z) = |\lam|e^{i\Th_j}\ov{z} + d_j$
otherwise). Thus, $A$ is the unique nonempty compact set satisfying
$$
A = \bigcup_{j=1}^m A_j,\ \ A_j = f_j(A).
$$
We further assume that the IFS $\{f_j\}_{j=1}^m$ is {\em invertible}
in the sense of Kameyama \cite{kam1}, that is, that there exists a continuous
map $q:\, A \to A$, such that $q|_{A_j}$ is the inverse of $f_j:\,A\to A_j$.
This means that there is ``compatibility'' of the IFS on the {\em overlap set}
$\Ok=\bigcup_{i\ne j} (A_i \cap A_j)$. Thus, we can consider the dynamical
system $(A,q)$, and $\Ok$ is the critical set of $q$ in the sense that $q$ is not a local homeomorphism 
precisely at $\Ok$.
We also assume that the overlap set is finite, which implies, by a recent
result of Bandt and Rao \cite{BR}, that the IFS is {\em critically non-recurrent}
(or simply non-recurrent). For an invertible IFS this means, by definition, that 
for every $z\in \Ok$, the limit set of $\{q^nz\}_{n\ge 1}$ does not contain $z$.
Being post-critically finite, of course, is a stronger property. Finally, we assume that $A$ is connected, so that
$\Ok$ is non-empty.

A rational map $p$ is called {\it semi-hyperbolic} if $p$ has no parabolic periodic
points and if all critical points in the Julia set $J$ are non-recurrent.

\begin{theorem} \label{thm-qs} Let $(A,q)$ be as above with finite non-empty overlap
set $\Ok$, and suppose that $\vphi: J\to A$ conjugates $p$ and $q$, 
where $p$ is a rational map and $J$ is its Julia set.
If $A$ is of bounded turning and if $p$ has no parabolic cycles, 
then $\vphi$ is quasisymmetric. 
\end{theorem}

Notice that under the assumptions of Theorem~\ref{thm-qs} $p$ is  semi-hyperbolic: Indeed, the critical points of $p$ on $J$ are precisely
those points where $p|_J$ is not a local homeomorphism, hence the conjugacy $\vphi$ takes critical points into critical points.

Intuitively, a homeomorphism is quasisymmetric if images of disks are ``roundish''
(see the discussion below for details). In order to prove quasisymmetry of
the conjugacy, we follow the strategy of  \cite{MS} in a non-hyperbolic
setting, as in \cite{PR}: The $\vphi$ image of a small disk
$B(x,r)$ centered at $x\in J$ is analyzed by first passing to the forward iterate
$p^{n}(B(x,r))$ for suitable $n$ (such that the diameter is large), second applying
$\vphi$ (noticing that on large scale, images of roundish sets stay roundish just by
continuity), and finally iterating backwards using $q$. That it is possible to pass
from small scale to large scale with bounded distortion in the critically non-recurrent
setting is due to Carleson, Jones and Yoccoz \cite{CJY}, whereas the same fact
in the linear setting is very simple. Now for the details.

Recall (see \cite{H}) that a homeomorphism
$f$ between metric spaces $(X,d_X)$ and $(Y,d_Y)$ is {\it quasisymmetric} (qs for short) if there is
a homeomorphism $\eta:[0,\infty)\to [0,\infty)$ such that for every $t>0$ and $a,b,x\in X$,
$d_X(x,a)\leq t\ d_X(x,b)$ implies $d_Y(f(x),f(a))\leq \eta(t)\ d_Y(f(x),f(b))$.
We will use the Polish notation $|a-b|$ instead of $d(a,b).$
The homeomorphism $f$ is called {\it weakly qs} if there is $H\geq1$ such that
$|a-x|\leq|b-x|$ implies $|f(a)-f(x)|\leq H|f(b)-f(x)|.$ By a result of Tukia and
V\"ais\"al\"a (\cite{TV}, \cite{H} Chapter 10), weak quasisymmetry implies quasisymmetry
if $X$ is doubling and connected (as subsets of $\R^2,$ our spaces are automatically doubling).
Recall that a space $X$ is of {\it bounded turning}, if
there exists a constant $L>0$
such that for every $z_1,z_2 \in X$ there is a curve $[z_1,z_2]\subset X$ satisfying
$$
\diam[z_1,z_2] \le L |z_1-z_2|.
$$
For $S\subset X$ and $x\in S$, denote the inradius of $S$ with respect to $x$ by
$\inrad(S,x) = \sup\{r\geq0: B(x,r)\subset S\} = \inf\{|x-y|:y\in X\setminus S\},$
where $B(x,r)$ is the closed
ball of radius $r$ centered at $x$. We say that a set $S\subset X$ is $C$-{\it roundish at x}
if $\diam S\leq C\ \inrad(S,x).$ Thus balls are 2-roundish (at their center).

\begin{lemma} \label{lem-qstest} Let $f:X\to Y$ be a homeomorphism and assume that $X$ is
of bounded turning. Then $f$ is qs if and only if the images of balls $B(x,r)$
under $f$ are uniformly roundish at $f(x)$.
\end{lemma}

\noindent
{\em Proof.} First assume that  $f$ is $qs.$ Let $y\in Y\setminus f(B(x,r))$ and $z\in f(B(x,r)),$
then  $|f^{-1}(z)-x|\leq r< |f^{-1}(y)-x|$ and hence $|z-f(x)|\leq H |y-f(x)|$
so that $\diam f(B(x,r))\leq 2H\ \inrad(f(B(x,r),f(x))$ by taking the supremum over $z$ and the infimum
over $y.$ For the converse, assume $\diam f(B(x,r))\leq M\ \inrad(f(B(x,r)),f(x))$ for all $x,r$.
If $|a-x|<|b-x|$, then
$|f(a)-f(x)|\leq \diam f(B(x,|a-x|)) \leq M\ \inrad(f(B(x,|a-x|),f(x))\leq M|f(b)-f(x)|.$
If $X$ is bounded turning, this implies weak quasisymmetry and hence quasisymmetry, as follows:
Given $a,b,x\in X$ with $|a-x|\leq|b-x|,$ let $[x,a]$ be a curve as in the definition of bounded turning.
Similar to the proof that weak qs implies qs (Theorem 10.19 in \cite{H}), it is easy to
construct a sequence of points $a_0,a_1,...,a_n$ on $[x,a]$ with $a_0=x, a_n=a$, such that
$|b-x|>|a_1-x|=|a_1-a_0|>|a_2-a_1|>\cdots>|a_n-a_{n-1}|$ and so that $n$ is bounded in terms of the
bounded turning constant $C$ and the doubling constant (a universal constant in $\R^2$) only.
Then
$$|f(a_i)-f(a_{i+1})|\leq M|f(a_i)-f(a_{i-1})|\leq M^{i+1} |f(b)-f(x)|,$$
and summation yields $|f(a)-f(x)|\leq n M^n |f(b)-f(x)|.$
\qed

\medskip

We will need the following version of results of Carleson, Jones and Yoccoz \cite{CJY}.

\begin{theorem} \label{thm-cjy} 
If $p$ is semi-hyperbolic, then
$J$ is of bounded turning. Furthermore, if $J$ is not the whole sphere, then
there are $\delta,r_0,c>0$ such that for every
$r<r_0$ and every $x\in J$ there is $n\geq 1$ such that
\begin{equation} \label{eq-cjy}
\Comp_x p^{-n} (B(p^n (x), c \delta)) \subset B(x,r)\subset \Comp_x p^{-n} (B(p^n (x), \delta))
\end{equation}
and such that the degree of $p^n$ on $\Comp_x p^{-n} (B(p^n (x), \delta))$ is bounded above
by the degree of $p$. In addition, we can choose $\delta< \delta_0$ for a given $\delta_0>0$.
\end{theorem}
Here $\Comp_x S$ denotes the connected component of $S$ containing $x.$
The paper \cite{CJY} deals with polynomials $p$. The bounded turning property
for rational functions can be found in \cite{Mih}, Corollary 2 (see the remark regarding
Corollary 2 in Section 4 of \cite{Mih}). For a proof of (\ref{eq-cjy}) in the rational case,
see the proof of Theorem A in \cite{PR} (using the notation of that proof, notice that 
$G(x,\delta,\operatorname{deg}(p),p)=\N$, use that the sequence 
$r_n=\diam \Comp_x p^{-n} (B(p^n (x), \delta))$ satisfies $r_n\leq \xi^n$ for some $\xi<1$
by \cite{PR}, Proposition 2.5, and notice that $r_{n+1}/r_n$ is bounded from below; it follows
that for every sufficiently small $r$ there is $n$ with $r_n$ comparable to $r$).

\noindent
\begin{lemma} \label{lem-diam} For a critically non-recurrent map $q$ there exist constants $C\ge 1$ and $\delta_1>0$ such that if 
$S\subset A$ is connected, then
\begin{equation}\label{eq-diam}C^{-1} (\frac1{|\lam|})^n\ \diam S 
\leq \diam q^n(S) \leq C (\frac1{|\lam|})^n\ \diam S,
\end{equation}
provided either $\diam q^n(S) < \delta_1$ or $(\frac1{|\lam|})^n\ \diam S < \delta_1$.
\end{lemma}

\noindent
{\em Proof.} 
Suppose that $T\subset A$ is connected. If $T\cap \Ok = \es$, then $T$ lies in one of the $A_j$'s, hence 
$\diam q(T) = \diam T/|\lam|$. If $T\cap \Ok \ne \es$ and 
$$\diam T < \delta_2:=\min_{x,y\in \Ok,\ x\neq y} |x-y|,$$
then $T\cap \Ok$ is a single point, say $z$, and 
\begin{equation} \label{eq-diam1}
\frac{1}{2|\lam|}\,\diam T \le \diam q(T) \le \frac{2}{|\lam|} \,\diam T.
\end{equation}
Indeed, the left-hand side follows by taking $x\in T$ such that $|z-x|\ge \half \diam T$, and the right-hand side follows by
writing, for $x_1,x_2\in T$: 
$$|q(x_1)-q(x_2)| \le |q(x_1)-q(z)| + |q(x_2) - q(z)|= |\lam|^{-1}(|x_1-z|+|z-x_2|).$$
By (\ref{eq-diam1}),
\begin{equation} \label{eq-diam2}
2^{-N_j}(\frac1{|\lam|})^j \diam T \le \diam q^j (T) \le 2^{N_j}(\frac1{|\lam|})^j \diam T,
\end{equation}
provided $\diam q^\ell(T) < \delta_2$ for all $\ell \le j$,
where $N_j=\{0 \le \ell \le j:\,q^\ell(T)\cap \Ok \ne \es\}$ is the number of times the iterates of $T$ ``hit'' the critical set.
By non-recurrence of $q$, there exists $\delta_3>0$ such that for every $z\in \Ok$ and every $\ell\ge 1$, $q^\ell(z)\not\in B(z,\delta_3)$.
If $z\in q^j(S) \cap q^{j+\ell}(S)$ for some $z\in \Ok$, then
$$
\delta_3 \le |q^\ell z- z| \le \diam q^{j+\ell}(S).
$$
Now suppose that
$$\diam q^n(S) < \delta_1:= \min\{\delta_2,\delta_3\}2^{-|\Ok|}.$$
By induction on $j$ (using (\ref{eq-diam2}) with $T = T_j = q^{n-j}(S)$), we conclude that
each $z\in \Ok$ can belong to $q^j(S)$ for at most one value of $j=0,\ldots,n$. Thus (\ref{eq-diam2}) yields (\ref{eq-diam}) with $C= 2^{|\Ok|}$.
The same argument applies when $(\frac1{|\lam|})^n\ \diam S < \delta_1$.
\qed

\medskip

The following lemma is an analog of (\ref{eq-cjy}) in the piecewise linear setting.
\noindent
\begin{lemma} \label{lem-roundish} Let $q$ be non-recurrent and assume that $A$ is of bounded turning, with the bounded turning constant
$L$.
Then 
\begin{equation}\label{eq-roundish} C\, |\lam|^n\,\diam B\ \ge \diam W \ge \inrad(W,x) \ge (CL)^{-1}\,|\lam|^n\, \inrad(B, q^n(x)),
\end{equation}
whenever $B\subset A$ is connected, $x\in q^{-n}(B)$ and $W=\Comp_x q^{-n}(B)$, provided that $\diam B < \delta_1$, where
$C$ and $\delta_1$ are from Lemma~\ref{lem-diam}.
\end{lemma}

\noindent
{\em Proof.} From Lemma \ref{lem-diam} we have
$$\diam B \geq \diam q^n(W)\ge C^{-1} (\frac1{|\lam|})^n \diam W.$$
If $w\in A$ is such that $|w-x|\leq \inrad(B,q^n(x)) |\lambda|^n / (CL)$, 
then there is a curve $\gamma$ joining $x$ and $w$ such that $\diam \gamma\leq L|x-w|$, and
Lemma \ref{lem-diam} yields $\diam q^n(\gamma)\leq \inrad(B,q^n(x)).$ It follows that $q^n(\gamma)\subset B$
and hence $\gamma\subset W$, in particular, $w\in W$. Thus
$$\inrad(W,x)\geq \inrad(B,q^n(x)) |\lambda|^n / (CL),$$ and the lemma follows. \qed

\medskip

\noindent
{\em Proof of Theorem~\ref{thm-qs}.} By Lemma \ref{lem-qstest}, it suffices to show
that images of balls are roundish. For balls of large radius this is true by continuity, so we can fix $x\in J$, $0<r<r_0,$ and let $n$
be obtained by Theorem \ref{thm-cjy}.
Because $\vphi$ is a conjugacy, we have
$$\vphi(\Comp_x p^{-n}(S)) = \Comp_{\vphi(x)} q^{-n}(\vphi(S))$$
for all sets $S$ and $x\in p^{-n}(S).$ Thus (\ref{eq-cjy}) gives
$$W_1=\Comp_{\vphi(x)} q^{-n} (\vphi(B(p^n (x), c\delta))) \subset \vphi(B(x,r))
\subset \Comp_{\vphi(x)} q^{-n} (\vphi(B(p^n (x),  \delta)))=W_2.$$
By continuity, $\vphi(B(p^n (x), c \delta)$ is uniformly roundish at $\vphi(p^n(x))$,
and its diameter is bounded away from zero and thus uniformly comparable to the
diameter of $\vphi(B(p^n (x), \delta)$. We can choose $\delta$ sufficiently small, so that the $\varphi$ image of any $\delta$ ball is contained
in a $\delta_1$ ball. Then
Lemma \ref{lem-roundish} shows that the diameters of $W_1$ and $W_2$ are comparable,
and that the inradius of $W_1$ is comparable to its diameter.
The theorem follows. \qed

\medskip

\noindent {\em Proof of Theorem~\ref{th-quad1}.} By Theorem~\ref{th-bt1}, $A_\lam$ is of
bounded turning, so the complement of $A_\lam$ is a John domain by \cite[Theorem 4.5]{NV}.
Theorem \ref{thm-qs} yields that the conjugacy $\vphi$ is quasisymmetric. (Note that parabolic cycles are excluded by the assumption that $J_c$ is a dendrite.)
Since quasisymmetric maps between boundaries of John domains extend to quasiconformal
homeomorphisms of the sphere, \cite[Theorem 4.2]{GNV}, it follows that $\vphi$ 
extends to a globally quasiconformal map. \qed

\section{Examples: rational maps}

We present two examples; most likely, one can find many more.

1. Let $A$ be the standard Sierpi\'nski gasket, the attractor of the IFS 
$$
\Bigl\{\half(z+1),\, \half(z + e^{2 \pi i/3}),\, \half (z + e^{4 \pi i/3})\Bigr\}\ \ \mbox{or}\ \ 
\Bigl\{\half(z+1),\, \half e^{2 \pi i/3}(z+ e^{2 \pi i/3}),\, \half  e^{4 \pi i/3} (z + e^{4 \pi i/3})\Bigr\}.
$$
The first IFS is the one commonly used, but the second IFS has the advantage of being invertible.
Let $q$ be the continuous map on $A$ whose inverse branches produce the second IFS.
It is known \cite{kam2} that $(A,q)$ is conjugate to $(J,p)$, where
$p(z) = z^2-16/(27z)$ and $J$ is its Julia set. (This claim is easy to verify: $p$ has three critical points, one of which is mapped to a fixed point, and another
two are mapped into a 2-cycle; this combinatorics agrees with that of the map $q$.) The rational map $p$ is semi-hyperbolic, and $A$ is obviously of bounded turning.
Thus, the conjugacy is quasisymmetric by Theorem~\ref{thm-qs}. 

2. Consider the ``hexagasket'' $H$ shown in Figure 4, the attractor of the
IFS $\{\frac{1}{3}(z-a_k) + a_k\}_{k=0}^5$ where $a_k = e^{\pi i k/3}$. It is, however, not invertible. Let $H =
\bigcup_{k=0}^5 H_k$ where $H_k$ is the piece of $H$ containing $a_j$. We observe that
$H$ can also be represented as a repelling invariant set 
of the piecewise linear function
$$
q(z) = \left\{ \begin{array}{llll} 3a_4 z - 2a_4, & z\in H_0, &
                                 3a_1 z -2a_2, & z\in H_1,\\
                                 3a_4z - 2a_0, & z\in H_2, &
                                 3a_1z - 2a_4, & z\in H_3,\\
                                 3a_4z - 2a_2, & z\in H_4, &
                                   3a_1z-2a_0, & z\in H_5. \end{array} \right.\ \
$$
The critical set consists of six points: $H_0\cap H_1 = \{b_1\}$,\ $b_1 = \frac{1}{\sqrt{3}}e^{\pi i/6}$,
$H_1\cap H_2 = \{b_2\}$, $b_2 = \frac{1}{\sqrt{3}}i$, etc., and one easily checks that $q$ is well-defined on them,
hence continuous. 
It is easy to check that $q$ maps the critical points as follows: $q(b_1)=q(b_4)=a_0$, $q(b_2)=q(b_5) = a_4$, and $q(b_3)=q(b_5)=a_2$. The critical values
form a 3-cycle: $a_0\mapsto a_4\mapsto a_2$; the remaining ``vertices'' $a_1, a_3, a_5$ are also mapped into this cycle.

Now consider the rational map $p(z) = z^4 + \frac{4}{9}e^{\pi i/3} z^{-2}$ whose Julia set $J$ is also shown in Figure 4. Note
that $p$ has six critical points $\zeta_j 2^{1/6}/3^{1/3}$ where $\zeta_j,\ j=1,\ldots,6,$ are the
sixth roots of $e^{\pi i/3}$. The critical values $a,b,c$ form a 3-cycle. Symmetry considerations (or a straightforward computation)  yield
that $p$ permutes the (straight line) rays to infinity through the critical values, and that $|p|$ is increasing on these rays.  Thus the
critical values are in the boundary of $G_\infty$, the unbounded component of $\C\setminus J$. 
A point $w$ in $G_\infty$ near a critical value, let's say near $b$, has six preimages under $p$. Four of them are in $G_\infty$ and two are in 
$G_0$, the component of $\C\setminus J$ containing $0$. 
One of the four is near $a$ ($q(a)=b$), and two are near the critical points that map to $b$. It follows that the critical points are in the
boundary of $G_\infty$. The two critical points that map to $b$ can contribute only one preimage each 
(the fourth preimage of $w$ in $G_\infty$ comes from a rotation of $b$ by $\pi/3$). So the other
two preimages of $w$ near the critical point are the two preimages in $G_0$, and it follows that the critical points are in the boundary of $G_0$.
This proves that the critical points are ``cut-points,'' and $P$ is a ``Misiurewicz-Sierpi\'nski'' map studied by Devaney et. al. \cite{Dev}, whose
Julia set $J(P)$ is a generalized Sierpi\'nski gasket. The conjugacy of $(H,q)$ to $(J,p)$ is well-defined on the critical orbits and extends by symbolic
dynamics to the entire set. By  Theorem~\ref{thm-qs}, it is quasisymmetric.

\begin{figure}[ht]
\begin{center}
$\begin{array}{cc}
\multicolumn{1}{l}{} &
        \multicolumn{1}{l}{} \\
\epsfxsize=2.in
\epsffile{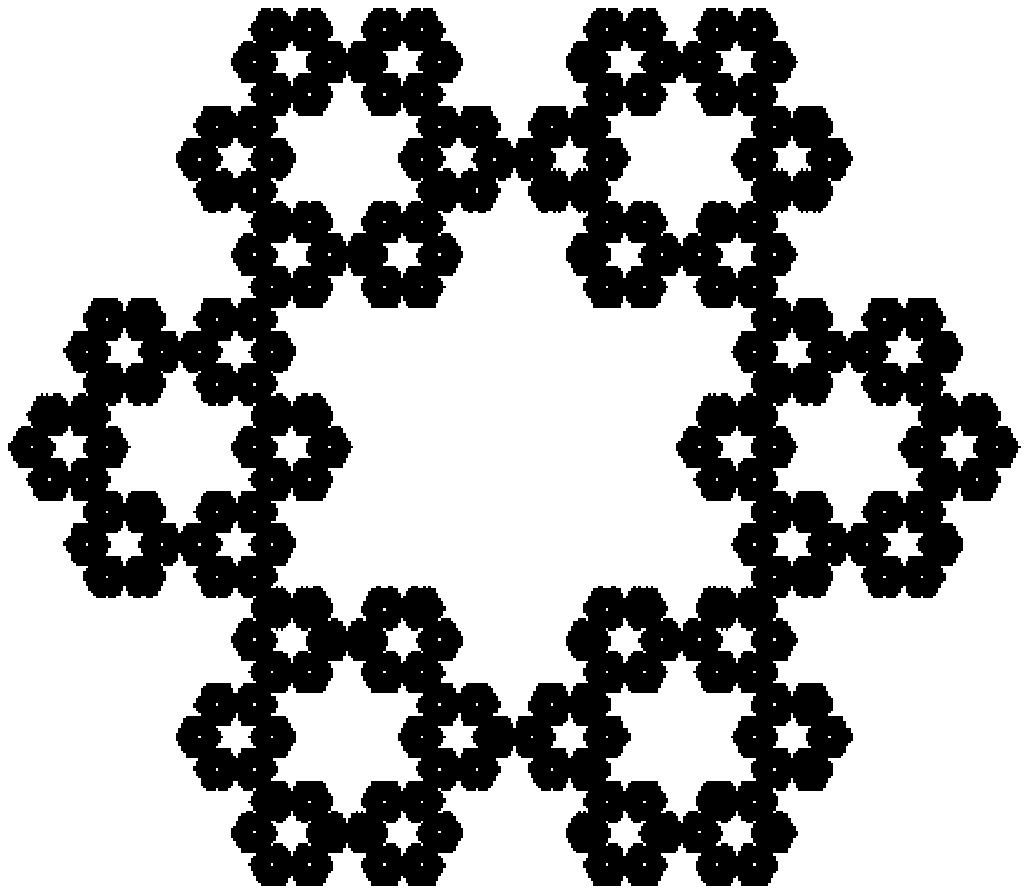} \ \ & \ \
        \epsfxsize=2.in
        \epsffile{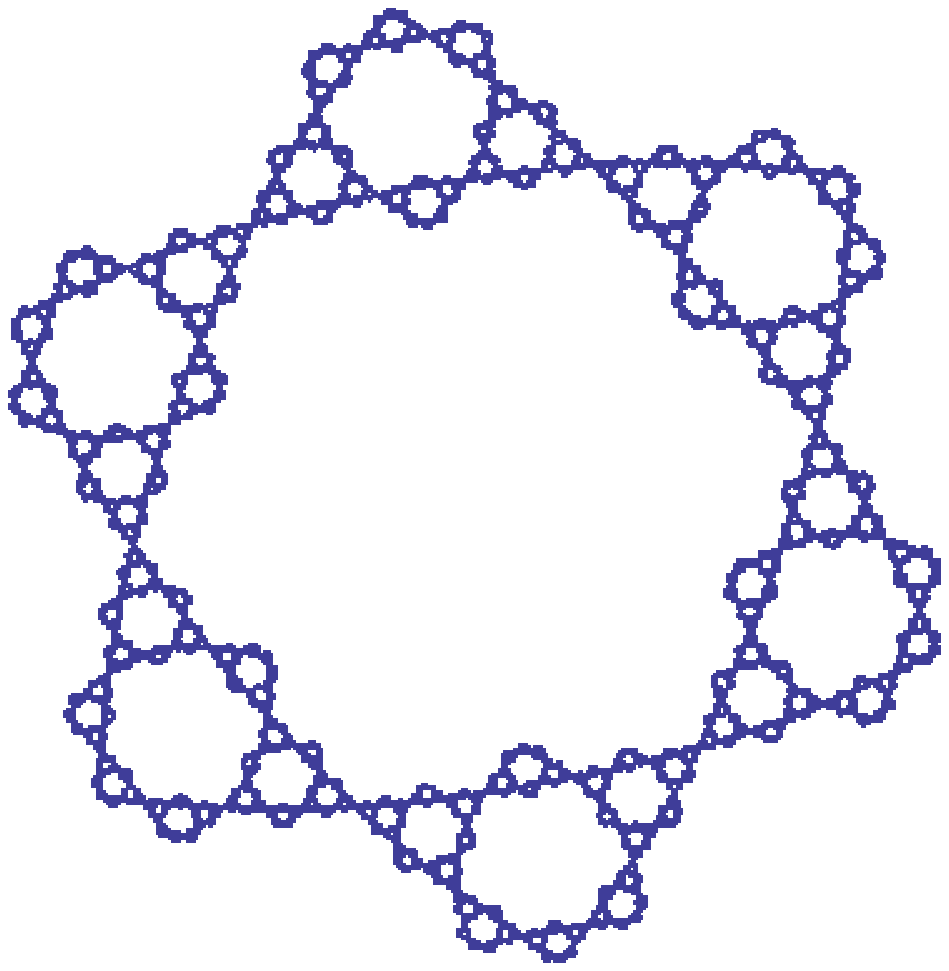} \\
\end{array}$
\end{center}
\caption{Self-similar set $H$ and the Julia set $J(P)$.}
\end{figure}

\section{Conjugacy}

Here we revisit the question of conjugacy between $(A_\lam,q_\lam)$ and $(J_c,p_c)$. Recall that $p_c(z) =
z^2+c$ for $c\in \C$, and $J_c$ is its Julia set.

Bandt \cite[Theorem 8.1]{bandt} claimed that if $A_\lam$ is a dendrite, then
there exists $c$ such that $(A_\lam,q_\lam)$ and $(J_c,p_c)$ are topologically conjugate.
However, his argument works only in the
post-critically finite case.
In fact, \cite{bandt} refers to a theorem of Kameyama who used
Thurston's Theorem on the topological characterization of
post-critically finite rational maps.
We review Kameyama's approach briefly, since it is not explained in \cite{bandt}.
We quote some
definitions from \cite{kam1}. 

\begin{defi}
{\em The pair $(K,\{F_i\}_{i=1}^d)$ is called a {\em self-similar system} if $K$ is
a compact set and there exists a continuous surjection
$\pi:\,\Sig_d\to K$ (the coding map) such that $F_i\pi(\om) = \pi(i\om)$ for all
$i$ and for all $\om \in \Sig_d$. 
The set $\Ok= \bigcup_{i\ne j} (K_i\cap K_j)$
is called the {\em overlap set}, where $K_i = F_i(K)$.
The system is {\em invertible} if
$F_i$ are injective and there exists a continuous map $g:\,K\to K$ such that $F^{-1}_i = g|_{K_i}$ for all $i$.
The self-similar system is {\em post-critically finite} if $\Ok$ is finite and
all its points are strictly preperiodic for $g$. We can consider $(K,g)$ as
a dynamical system with the critical set $\Ok$.}
\end{defi}

We have already seen invertible IFS in a special case, at the beginning of the last section. Note that
``self-similarity'' in the definition above is very different from self-similarity used in this paper; the
former is just a way to define the IFS, whereas the latter requires that the maps are similitudes.
Note also that Kameyama uses the term ``connecting set'' instead of ``overlap set'' used here.

Kameyama \cite[Th.B]{kam1} (see also \cite[Th.6.1]{kam2}) proved that if
$(K,\{F_i\}_{i=1}^d)$ is an invertible post-critically finite self-similar
system, with $K$ connected and simply connected, then $(K,g^n)$ is
topologically conjugate to $(J_P,P)$ for some $n$, where $P$ is a polynomial
and $J_P$ is its Julia set. Moreover, $n=1$ if $g$ preserves the cyclic
order of ``arms'' at the vertices of the Hubbard tree, defined as the connected hull of the critical orbit.
The vertices are defined as the branch points of the tree intersected with $K_i$, $i\le d$.

In our case, for $\lam \in \Tk$,
$(A_\lam,\{F_0,F_1\circ s\})$ is an invertible self-similar system, with $g=q_\lam$, see (\ref{eq-ql}).
The dendrite $A_\lam$ is connected and simply connected, so in the post-critically finite case Kameyama's Theorem
applies. We have $n=1$ because the action of $q_\lam$ on the
subsets $A_i$ is linear, so the cyclic order at branch points is preserved. Thus, $(A_\lam,q_\lam)$ is
conjugate to $(J_P,P)$ for some polynomial, which has to be quadratic and can be chosen in the form
$p_c(z) = z^2+c$.

\subsection{Symbolic dynamics of $q_\lam$.}
Next we present an alternative approach to the question of conjugacy, based on
symbolic dynamics. We are grateful to Henk Bruin for many helpful comments and suggestions which are used in
this section.

Let $\Sig_2 = \{0,1\}^\N$; denote by $\sig$ the left shift. We fix $\lam \in \Tk$ and let 
$$
\pi_\lam(\ba) = \sum_{n=0}^\infty a_n \lam^n,\ \ \ \mbox{for}\ \ba = a_0a_1\ldots \in \Sig_2.
$$
We have $A_\lam = \pi_\lam(\Sig_2)$. The sequence $\ba$ is the symbolic ``address'' of $z=\pi(\ba)\in A_\lam$,
corresponding to the IFS $\{F_0,F_1\}$. In order to determine the $q_\lam$-itinerary, we need to ``rewrite'' this address
in terms of the IFS $\{F_0,F_1\circ s\}$. Let $s$ be the ``flip'' map on $\Sig_2$ which switches $0$'s and
$1$'s.
From the definition (\ref{eq-ql}) of $q_\lam$ it follows that
\begin{equation} \label{eq-qact}
q_\lam(\pi(\ba)) = \left\{ \begin{array}{ll} \pi_\lam(\sig\ba), & \mbox{if} \
                        a_0 =0; \\
                    \pi_\lam(s(\sig\ba))=\pi_\lam(\sig(s(\ba))), & \mbox{if} \ a_0=1.
\end{array} \right.
\end{equation}
That is, $q_\lam$ acts on the address as a left shift or  a shift and a flip, depending on the digit $a_0$.

We assign the symbols $0$, $1$, and $\star$ to $A_0\setminus \{o_\lam\}$, $A_1 \setminus \{o_\lam\}$, and
$\{o_\lam\}$ respectively, and define itineraries and the kneading sequence $\nu_\lam$, in the same way as for
$p_c$ on a connected, locally connected Julia set. By (\ref{eq-qact}), if the orbit of $z$ does not hit $o_\lam$, then
the {\em itinerary} $e(z) = e_0e_1 e_2\ldots \in \Sig_2$ of $z = \sum_{n=0}^\infty a_n \lam^n$ satisfies
$$
e_0=a_0\ \ \ \mbox{and for $n\ge 1$:}\ \ e_n = \left\{ \begin{array}{ll} 1 & \mbox{if}\ a_n \ne a_{n-1}, \\
                                                                         0 & \mbox{if}\ a_n = a_{n-1}
\end{array} \right.
$$
Conversely, $a_n =0$ or $1$ according to whether $e_0\ldots e_n$ contains an even number of $1$'s. 
The {\em kneading sequence} is basically the itinerary of the critical value. However, it is
customary to start it with a $1$, so it is defined by $\nu = \nu_\lam = e(q_\lam(o_\lam))$ or $s(e(q_\lam(o_\lam)))$
according to whether $q_\lam(o_\lam)$  is in $A_1$ or $A_0$. Explicitly:
\begin{equation} \label{eq-knead}
\nu = \nu_1\nu_2\ldots, \mbox{where}\ \nu_n = \left\{ \begin{array}{ll} e_{n-1}(q_\lam(o_\lam)), & \mbox{if} \
                                              q_\lam(o_\lam)\in A_1; \\
                                            1-e_{n-1}(q_\lam(o_\lam)), & \mbox{if} \
                                              q_\lam(o_\lam)\in A_0. \end{array} \right.
\end{equation}
As was already mentioned, by a result of Bandt and Rao \cite{BR}, in our case ($\lam \in \Tk$, i.e.\
$A_\lam$ is a dendrite)
the system $(A_\lam,q_\lam)$ is critically non-recurrent, hence the orbit of $o_\lam$ does not return
to $o_\lam$ and the kneading sequence is in $\Sig_2$.

As shown in \cite{solo}, $\lam\in \Tk$ if and only if there is a unique powers series $g_\lam(z)=
1+\sum_{n=1}^\infty b_n z^n$, with $b_n \in \{-1,0,1\}$, satisfying $g_\lam(\lam) = 0$, and moreover, all
$b_n$ are in $\{-1,1\}$. This power series is obtained from the address of the critical point:
$$
o_\lam = \sum_{n=0}^\infty a_n \lam^n = \sum_{n=0}^\infty (1-a_n) \lam^n, \ \mbox{with}\ a_0=1\ \ \Longleftrightarrow\ \ 
b_n = 2a_n -1,\ n\ge 1.
$$  
In view of the above, the sequence $b_n$ is expressed in terms of the kneading sequence as follows:
\begin{equation} \label{eq-knead2}
b_n = \left\{ \begin{array}{ll} (-1)^{\#\{1 \le i \le n:\ \nu_i = 1\}} , & \mbox{if} \
                                              q_\lam(o_\lam)\in A_1; \\
                                 (-1)^{\#\{1 \le i \le n:\ \nu_i = 0\}} , & \mbox{if} \
                                              q_\lam(o_\lam)\in A_0. \end{array} \right.
\end{equation}
Once we have the power series $g_\lam$, the parameter $\lam$ is determined (essentially) uniquely, since we know by \cite{BBBP} that a power series with coefficients $\pm 1$
can have at most one
zero in $\{z=x+iy:\ y>0,\ |z|^2 \le \half\}$. The two power series obtained from (\ref{eq-knead2}) will be related by $\wtil{g}(z) = g(-z)$, so will correspond to the choice
of $\lam$ or $-\lam$. Since the maps $q_\lam$ and $q_{-\lam}$ are conjugate (the conjugacy is $\Phi(z) = -z+\frac{1}{1-\lam^2}$ from $A_\lam$ to $A_{-\lam}$), 
we can always assume, without loss of generality that $b_1=1$ and $q_\lam(o_\lam)\in A_1$.

\subsection{Connection with quadratic dynamics.}
Suppose that $J_c$, the Julia set of $p_c(z) = z^2+c$, is a dendrite. Then there is a standard way to write
$J_c = J_0\cup J_1 \cup \{0\}$ (a disjoint union),   using external rays, so that $p_c(0)\in J_1$.
Then the kneading sequence of $p_c$ is defined as the itinerary of the critical value $c$.

\begin{prop} \label{prop-conj} Let $\lam \in \Tk$ and $c\in \C$ is such that $J_c$ is a dendrite,
and the kneading sequences of $q_\lam$ and $p_c$ are identical. Then $(A_\lam,q_\lam)$ and $(J_c,p_c)$ are
conjugate.
\end{prop}

\noindent {\em Proof.} Both systems
have the same ``symbolic model'': the sequence space $\Sig_2$ factored by the equivalence relation:
$\ba \sim \ba'$ iff there exists $n\ge 0$ such that $a_j = a'_j$ for $j< n$ and $\sig^n\ba$, $\sig^n\ba'$
are the two addresses of the critical point $o_\lam$. We map a point in $A_\lam$ with a given $q$-itinerary
into the point in $J_c$ with the same itinerary. We will show that this map, let's denote it $\phi$, is well-defined.
Indeed, it is easy to see that $p_c$ is  critically non-recurrent since the kneading sequence has an initial block which never occurs again.
We assumed that $J_c$ is a dendrite, in particular, that it is locally connected. 
This implies that the conformal map from $\C^*\setminus \D$ to $\C^*\setminus J_c$ (where $\C^*$ denotes the Riemann sphere), conjugating $z^2$ to $p_c$, which is used to define
symbolic dynamics on $J_c$, extends continuously to the
unit circle, and hence the ``cylinder sets'' in $J_c$, defined as sets of points with a common itinerary of length $n$, shrink to a singleton, as $n\to \infty$.
This proves that $\phi$ is well-defined and continuous. It is straightforward to check that it is a conjugacy. \qed

\medskip

In view of the proposition, the following strategy seems natural: given $\lam \in \Tk$, determine the
kneading sequence $\nu_\lam$ of $q_\lam$, and try to find the quadratic map with a dendrite $J_c$, with the same kneading sequence.

There is a connection between the kneading sequences of quadratic maps and those for the doubling map on the circle.
Suppose that $J_c$ is locally connected and the critical value $c$ is a landing point for an external ray with angle $\theta$. Then all dynamical
external rays land at points of $J_c$, see \cite[Theorem 18.3]{Milnor}. The rays with angles $\theta/2$ and $(\theta+1)/2$ land at the critical point $0$
and divide
$J_c\setminus \{0\}$ into the ``pieces'' $J_0$ and $J_1$, with $c\in J_1$. Now, the itinerary of $c$ under $p_c$ with respect to this partition coincides
with the itinerary of $\Th$ under the angle-doubling
map on the circle $S^1$ with respect to the partition $\{(\frac{\Th}{2},\frac{\Th+1}{2}), (\frac{\Th+1}{2},
\frac{\Th}{2})\}$ of  $S^1$, with the the arc $(\frac{\Th}{2},\frac{\Th+1}{2})$ corresponding to $1$, and
the complementary one to $0$. We call this the kneading sequence of the angle-doubling corresponding to $\Th$.

A sequence in $\Sig_2$ is called {\em admissible} if it is a kneading sequence for the angle-doubling map for some $\Th\in S^1$.
It is well-known that not every sequence is admissible. There are several criteria for admissibility,
see the book in progress by Bruin and Schleicher \cite{BS} (see also \cite{BS2,Kaffl}). We do not discuss them
here in detail. Sometimes, the following simple-minded approach works.

\begin{lemma} \label{lem-admis} Let $\bw\in \Sig_2$, with $w_1=1$, and
suppose that 
\begin{equation} \label{eq-admis}
\bw >_{lex}\sig^k\bw,\ \ \mbox{ for all}\  k\ge 1, 
\end{equation}
where $>_{lex}$
is the lexicographical order on $\Sig_2$. Then $\bw$ is admissible.
\end{lemma}

\noindent {\em Proof.} This is proved in \cite{BS}, but we present a short direct proof.

Consider $\Th = \sum_{n=1}^\infty (1-w_n) 2^{-n}$ and its kneading sequence for the  angle-doubling, as defined above. 
We claim that this is precisely the sequence $\bw$. Indeed, the iterates
of $\Th$ under the doubling map have binary expansions $(1-w_{k})(1-w_{k+1})\ldots$ for $k\ge 1$. If $w_k=1$, 
then the corresponding number $.0(1-w_{k+1})(1-w_{k+2})\ldots$ is greater than $\Th/2 = .0(1-w_1)(1-w_2)\ldots$ by (\ref{eq-admis}) and
$\le .1$ (in binary),
and so will be assigned the digit 1 in the itinerary. If $w_k=0$, then the
corresponding number $.1(1-w_{k+1})(1-w_{k+2})\ldots$ is greater than $\frac{\Th+1}{2} =
.1(1-w_1)(1-w_2)\ldots$ and $\le .\ov{1}$ (in binary), and so will be assigned the digit 0 in the itinerary. This concludes the proof of the claim, and
of the lemma. \qed 

\medskip

We summarize the discussion in the following statement.

\begin{prop} \label{prop-flaky} Suppose that $\lam\in \Tk$ is such that 

(i) the kneading sequence $\nu_\lam$ of $q_\lam$ is admissible, and corresponds to 
an angle $\Th$;

(ii) 
there is a quadratic Julia set $J_c$ which is locally connected, with an external angle $\Th$.

\noindent
Then $(A_\lam,q_\lam)$ and $(J_c,p_c)$ are conjugate.
\end{prop}

This is almost immediate from Proposition~\ref{prop-conj}; we only need to note that the kneading sequence is non-recurrent, hence $p_c$ is critically
non-recurrent, and then $J_c$ can only be a dendrite.

\begin{remark} \label{rem-land}
{\em It is folklore, although we do not know a reference, that 
in the above proposition (i) implies (ii).} 
\end{remark}

\begin{remark} \label{rem-Mis}
{\em If $q_\lam$ is post-critically finite, then $\nu_\lam$ is strictly preperiodic. It is admissible, since we
have conjugacy by Kameyama's Theorem, and
the corresponding angle $\Th$ is rational, with an even denominator. In this case
it is well-known \cite{DH} that the external parameter ray of angle $\Th$ lands at a {\em Misiurewicz point} $c$, for which $J_c$ is a dendrite.}
\end{remark}

\subsection{The work of Bruin and Schleicher}
Our systems $(A_\lam,q_\lam)$ satisfy all the axioms of an abstract Julia set of
Bruin and Schleicher \cite[Def.\,14.4]{BS}.
The connected hull of the orbit of $o_\lam$ in $A_\lam$ is the (abstract) Hubbard tree or Hubbard dendrite
according to whether the orbit is finite or infinite. In the post-critically finite case, we get a Hubbard
tree. It is proved in \cite{BS} (see also \cite[Prop.\,2.11]{Kaffl}) that a Hubbard tree corresponds to a
quadratic polynomial if and only if the ``arms'' are cyclically permuted at periodic branch points.
As already mentioned at the beginning of this section, $q_\lam$ has this property, and this yields another
proof of the conjugacy in the post-critically finite case.

Using the techniques of \cite{BS}, one can derive various properties of the dendrites from their
kneading sequence. In particular, in all our examples (see below),
but not for general kneading sequences, the following holds [personal communication from Henk Bruin]:
\begin{itemize}
\item all the $\alpha$-fixed point of $A_\lam$ have three arms, and there are no other periodic branch points;
\item the critical value $q(o_\lam)$ is an endpoint of $A_\lam$.
\end{itemize}
We do not know if this holds in general. The following questions are also open:
\begin{itemize}
\item Is it true that $\Tk\subset \Mk$ is a Cantor set?
\item Characterize those $c$ in the quadratic Mandelbrot set for which there exists a corresponding
$\lam \in \Tk$.
\end{itemize}

\subsection{Examples}
1. Let $
\lam\approx .3668760 + .5202594i,
$
be the (only) zero of the power series with coefficients
$\bB=+(-+++--)$ in the upper half-plane of modulus less than $1/\sqrt{2}$; here $+$ and $-$ correspond to $\pm 1$'s and
parentheses indicate a period.
By (\ref{eq-knead}), we get the kneading sequence
$\nu_\lam=1100100\ldots = 1\ov{100}$. The condition (\ref{eq-admis})
holds since $11$ never occurs in $\bw$ except at the beginning. Note that
the angle with the corresponding (0-1 flipped) binary expansion is $\Th=(\frac{1}{8}+\frac{1}{16})
\frac{8}{7}=\frac{3}{14}$. We conclude that $(A_\lam,q_\lam)$ is conjugate to
$(J_c,p_c)$ where $c$ corresponds to the point on the quadratic Mandelbrot set
with the external angle $\frac{3}{14}$. See Figure 2 for the picture of these sets.
In Figure 3 we show
the dendrite for $\lam$ with the power series $+--(+)$, with the
kneading sequence $\nu_\lam=101\ov{0}$, and the Julia set
corresponding to the external angle $\frac{3}{8}$.

2. Denote 
\begin{equation} \label{def-T0}
\Tk_0 := \{\lam:\ |\lam| <1/\sqrt{2}, \ \exists\,f\in \Bk_0,\ f(\lam)=0\},
\end{equation}
where $\Bk_0$ is the set of power series whose coefficients agree with $\bB=+(-+++--)$ up to degree 14, and after degree
14 consist of arbitrary concatenations of blocks $+++$, $++++$, $---$, $----$, with the blocks of $-$'s followed
by the blocks of $+$'s and vice versa.
In \cite{solo} it is proved that $\Tk_0\subset \Tk$, that is, $A_\lam$ with $\lam\in \Tk_0$ are all dendrites.
The family $\Tk_0$ is clearly uncountable (of positive Hausdorff dimension).

More examples can be found numerically.
We found (without rigorous proof, but with a high degree of confidence)
that zeros of power series ranging from
$\gam_1\approx 0.28093+ 0.59621 i$, a zero of power series
expressed by $+(-++--)$, to
$\gam_2\approx 0.486036+ 0.453914i$, a zero of power series
expressed by $+(-++++++----)$, are all in $\Tk$, as well as
a great many in between (they can be surmised from Figure 9(a)
in \cite{bandt}).
In all our examples (including $\lam\in \Tk_0$) the kneading sequence is admissible.
This follows from Lemma~\ref{lem-admis}
since the kneading sequence will start with $11$, and
$11$ will never occur again ($11$ in a kneading sequence means two sign changes on a row, that is, $+-+$ or $-+-$).
It may be that (\ref{eq-admis}) is always satisfied for $\nu_\lam$, with $\lam \in \Tk$, but we do not know
how to prove this.

3. It seems that most dendrite Julia sets $J_c$ do not correspond to any $A_\lam$ with $\lam\in \Tk$. For example, consider, perhaps, the best known dendrite Julia
set, with $c=i$. The critical value $i$ gets into the periodic 2-orbit $\{i-1,-i\}$ after one iteration, so we have the kneading sequence $\nu = 1 \ov{10}$ (the
external angle is $\frac{1}{6}$). If there
is $\lam\in \Tk$ with such a kneading sequence, then $\lam$ is a zero of a power series with the coefficients $\bB = +(-++-)$. This is a rational function of
degree four, whose zeros are the zeros of $1-z+z^2+z^3-2z^4= (z-1)(z+1)(1-z+2z^2)$. The zeros inside the unit disk are $\frac{1\pm \sqrt{7}i}{4}$, which have
modulus $2^{-1/2}$. In fact, the corresponding $A_\lam$ is known as ``tame twindragon,'' which tiles the plane by translations and has nonempty interior, 
so it is definitely not a dendrite. Another example is given in \cite[Remark 8.2]{bandt}: for the Julia set with external angle $\frac{1}{8}$, the only candidate
$\lam$ is $\frac{1\pm i}{2}$, which gives the ``twindragon'' $A_\lam$, again with nonempty interior. The
kneading sequence is $\nu=111\ov{0}$. Other external angles, for instance, $\frac{1}{10}$ or $\frac{3}{10}$,
yield unremarkable $\lam$ of modulus greater than $2^{-1/2}$, so again not in $\Tk$.

\section{Bounded turning property II}

Here we prove the bounded turning property for a large class of self-similar sets, which include p.c.f.\ $A_\lam$, not necessarily of
dendrite type. 

We consider iterated function systems of the
form $\{f_j(z)=\lam z+d_j\}_{j=1}^m$ for $m\ge 2$ and $\lam \in\D$.
The attractor $E$ has an explicit representation as the set of sums of
power series in
$\lam$ with coefficients in $\Dk = \{d_1,\ldots,d_m\}\subset \C$:
\begin{equation} \label{eq-attr}
E = \bigcup_{j=1}^m E_j = \Bigl\{\sum_{n=0}^\infty a_n \lam^n: a_n \in \Dk
\Bigr\}.
\end{equation}
We have $E = \bigcup_{j=1}^m E_j$, where $E_j=f_j(E)$ are the ``pieces'' of
of the attractor. Our main assumption is that $E$ is connected and
the ``overlap set''
$$\Ok:= \bigcup_{i\ne j}(E_i \cap E_j)$$ is finite. (Note, however, that here we do not assume the IFS to be invertible, as in 
Sections 2 and 3.) Unfortunately, we do not
know if this already implies the bounded turning property, so we have to
impose an additional technical assumption. We will show, however, that this
assumption is satisfied in many cases.
Let
$$
\Bk= \Bigl\{f(z)= \sum_{n=0}^\infty c_n z^n:\ c_n \in \Dk - \Dk,\ c_0\ne 0\Bigr\}.
$$
It is immediate that
$$
E_j \cap E_k \ne \es \ \Longleftrightarrow
\ \exists\,f\in \Bk,\ \ c_0=a_j-a_k,\ f(\lam) = 0.
$$
Let
\begin{equation} \label{eq-Fk}
\Fk =\Bigl\{f\in \Bk:\ f(\lam)=0\}.
\end{equation}
Since $E$ is connected, we have that $\Fk \ne \es$.

There is a natural projection $\pi:\, \Dk^\N\to A$ defined by
$\pi(\ba) = \sum_{n=0}^\infty a_n \lam^n$. The elements of $\pi^{-1}(\{x\})$
are called addresses of $x$.

Now we introduce our technical assumption. For two power series
$f,g$ denote by $|f\wedge g|$ the degree of their maximal common initial part,
i.e.\ $|f\wedge g|= \min\{k\ge 0:\ f^{(k+1)}(0) \ne g^{(k+1)}(0)\}$.

\begin{theorem} \label{th-bt}
Assume that the self-similar set $E$ of the form (\ref{eq-attr})
is connected and has a finite overlap set. In addition, suppose that there
exists $C_1>0$ such that for all $n\in \N$ sufficiently large,
\begin{equation} \label{eq-tech}
(f\in \Fk,\ g\in \Bk,\ |f\wedge g| = n)\ \Rightarrow\ |g(\lam)| \ge
C_1|\lam|^n.
\end{equation}
Then $E$ has bounded turning.
\end{theorem}

We say that $E$ is p.c.f.\ (post-critically finite) if the overlap set
is finite and
every sequence in $\pi^{-1}(\Ok)$ is eventually periodic.
This definition is consistent with the definition of p.c.f.\ used earlier, but it is more general, since here we do not
assume that the IFS is invertible.

\begin{corollary} \label{cor-pcf}
Assume that the self-similar set $E$ of the form (\ref{eq-attr})
is p.c.f.\ and connected. Then $E$ has bounded turning.
\end{corollary}

In addition to the p.c.f.\ case, Theorem~\ref{th-bt} covers the uncountable set $A_\lam,\ \lam\in \Tk_0$
defined in (\ref{def-T0}), since (\ref{eq-tech}) for those sets immediately follows from the estimates in 
\cite[Prop.\,4.1]{solo}.

\begin{corollary} \label{cor-new}
Let $\lam\in \Tk$ be such that the conditions of Proposition~\ref{prop-flaky} are fulfilled, with some $c=c(\lam)$,
and the self-similar set
$A_\lam$ satisfies (\ref{eq-tech}). Then $(A_\lam,q_\lam)$ and $(J_c,p_c)$ are quasisymmetrically conjugate.
In particular, this holds for all $\lam \in \Tk_0$.
\end{corollary}

This is just a combination of Theorem~\ref{thm-qs}, Proposition~\ref{prop-flaky}, and Theorem~\ref{th-bt}.
The claim about $\Tk_0$ follows from Remark~\ref{rem-land}.

\begin{remark} \label{rem-john}
{\em By \cite[Th.\,4.5]{NV}, the bounded turning property of a compact set $E$ implies that 
the components of the complement of $E$ are John domains.
In the paper \cite{Ai} it is shown that there exist self-similar sets in the plane (not of the kind considered in this paper) which satisfy
the Open Set Condition, but the complement does not have the John property. On the other hand, \cite{Ai} 
obtained some geometric conditions (which do not hold in our case) for the complement of a self-similar
set to be a uniformly John domain. The picture of a self-similar set in \cite[Fig.\,4]{BR} suggests that it does not have the bounded turning property, even though 
the overlap set is finite. The IFS, however, contains both orientation-preserving and orientation-reversing maps, so it is also of a different kind from those 
considered in our paper. } 
\end{remark}

For the proof of Theorem~\ref{th-bt}, we need a lemma, which
does not rely on (\ref{eq-tech}).

\begin{lemma} \label{lem-vspom}
Suppose that the self-similar set $E$ from (\ref{eq-attr})
has a finite overlap set $\Ok$. Then the set $\Fk$ from (\ref{eq-Fk}) is
finite, and moreover, for each power series $f\in \Fk$, the coefficients
$c_n = f^{(n)}(0)/n!$ have a unique representation $c_n = d_{j(n)} - d_{k(n)}$
for $d_{j(n)}, d_{k(n)}\in \Dk$, for all $n$ sufficiently large.
\end{lemma}

{\em Proof of the lemma.}
Bandt and Rao \cite{BR} proved that an IFS in the plane
whose attractor $A$ is connected and the overlap set is finite, satisfies the
Open Set Condition. It is well-known that this implies that no point in $A$
can have infinitely many distinct addresses.

If $f(z) = \sum_{n=0}^\infty c_n z^n\in \Bk$ is
such that $f(\lam)=0$ (that is, $f\in \Fk$) and $c_n = d_{j(n)} - d_{k(n)}$
for $d_{j(n)}, d_{k(n)}\in \Dk$, then
$$
\sum_{n=0}^\infty d_{j(n)} \lam^n = \sum_{n=0}^\infty d_{k(n)} \lam^n \in
A_{j(0)}\cap A_{k(0)} \subset \Ok.
$$
If a coefficient $c_n$ has non-unique representation as a difference of
elements in $\Dk$, say $c_n = d_{j(n)} - d_{k(n)}=d'_{j(n)} - d'_{k(n)}$, with $d_{j(n)}\ne d'_{j(n)}$, then we get a pair of distinct points in $\Ok$, with a
difference of $(d_{j(n)}-d'_{j(n)})\lam^n$.
If there are infinitely many such $n$, we get that $\Ok$
is infinite, a contradiction. This proves the second claim of the lemma.

In order to see that $\Fk$ is finite, note that every $f\in \Fk$ yields
at least one point in $\Ok$. Different functions in $\Fk$ yield either
distinct points in $\Ok$ or one point in $\Ok$ with multiple addresses.
It follows that $\Fk$ must be finite.
 \qed

\medskip

{\em Proof of Theorem~\ref{th-bt}.}
Consider $z_1\ne z_2\in E$, with the goal of proving (\ref{eq-bt}). We can
assume that $z_1$ and $z_2$ do not have addresses starting with the same
symbol, that is, there is no $j\le m$ such that $z_1, z_2 \not\in E_j$. Indeed, otherwise
we could pass to $f_j^{-1}(z_i)$, $i=1,2$, keeping in mind that
$f_j(z) = \lam z+ d_j$ is a similitude.
This process has to stop since $f_j^{-1}$ is an expansion on $E_j$.

We will further
assume that $z_1$ and $z_2$ are in $E_1$ and $E_2$ respectively, and
both are very close to some $x\in E_1\cap E_2$, since otherwise (\ref{eq-bt})
trivially holds with some constant $L$.

\smallskip

{\bf Claim.} {\em For every $K>0$ there exists
$\delta_1(K)>0$ such that every address of every
$z\in B(x,\delta_1)$ coincides with an address of $x$ in at least $K$ digits.
}

\smallskip

Indeed, if the claim is false, we have a sequence of points $z_n \to x$ and
addresses $\bB^{(n)} \in \Dk^\N$ such that $\pi(\bB^{(n)}) = z_n$
and no address of $x$ agrees with any of $\bB^{(n)}$ in the first $K$ digits.
Passing to a subsequence we can assume that $\bB^{(n)} \to \ba$, whence
$\pi(\ba) = x$, so $\ba$ is an address of $x$ and $\bB^{(n)}$ must eventually
agree with $\ba$ in the first $K$ digits. This contradiction proves the claim.

\smallskip

Recall that $x$ has finitely many addresses, with at least two of them
starting with 1,2 respectively. Let $\ell\in \Nat$ be so large that

\smallskip

(i) no two
distinct addresses of $x$ agree in all of the first $\ell$ digits;

\smallskip

(ii) for every $f\in \Fk$ and every $n\ge \ell$, the Taylor coefficient
$c_n = f^{(n)}(0)/n!$ has a unique representation as a difference of
elements in $\Dk$ (which is possible by Lemma~\ref{lem-vspom}).

\smallskip

(iii) the condition (\ref{eq-tech}) holds for $n\ge \ell$.

\smallskip

Let $\delta = \delta_1(\ell)$.
If either $z_1$ or $z_2$ are not in $B(x,\delta)$, then (\ref{eq-bt}) holds
with some constant $L$.
Now suppose that $z_1, z_2 \in B(x,\delta)$.
By the claim, there are addresses $\bB^{(i)}$ of $z_i$,
for $i=1,2$, such that $\bB^{(i)}$ agrees with
$\ba^{(i)}$ in at least $\ell$ first digits, where $\ba^{(i)}$ is
an address of $x$ starting with $i$. (Note that $z_i \ne x$, since otherwise
both points are in $A_j$ where $j\ne i$.)
We have that
$$
f(z) := \sum_{n=0}^\infty (a^{(1)}_n - a^{(2)}_n)z^n \in \Fk,\ \ \
g(z) := \sum_{n=0}^\infty (b^{(1)}_n - b^{(2)}_n)z^n \in  \Bk.
$$
Let $k_i$, $i=1,2$, be the length of the common initial segment of $\ba^{(i)}$ and $\bB^{(i)}$.
Assume that $k_1 \le k_2$ without loss of generality.
Observe that
$$
c_{k_1+1} =
a^{(1)}_{k_1+1} - a^{(2)}_{k_1+1} \ne b^{(1)}_{k_1+1} - b^{(2)}_{k_1+1},
$$
since $a^{(1)}_{k_1+1}\ne b^{(1)}_{k_1+1}$, and $c_{k_1+1}$ has
a unique representation as a difference of
elements in $\Dk$ by the assumption (ii). It follows that $|g\wedge f| = k_1$. Now, by (\ref{eq-tech}),
\begin{equation} \label{eq1}
|z_1-z_2| = |\pi(\bB^{(1)}) - \pi(\bB^{(2)})| =
|g(\lam)| \ge C_1|\lam|^{k_1} \ge (C_1/2) (|\lam|^{k_1} + |\lam|^{k_2}).
\end{equation}
Let $u^{(i)}$, $|u^{(i)}|=k_i$,
be the common initial part of $\bB^{(i)}$ and $\ba^{(i)}$.
We know that $A_{u^{(i)}}$ is connected
(being a similar copy of $A$), hence
$$
\diam[z_i,x] \le \diam A_{u^{(i)}}=|\lam|^{k_i} \diam A.
$$
Combining this with (\ref{eq1}) yields
$$
|z_1-z_2| \ge (C_1/2) (\diam A)^{-1} (\diam[z_1,x] + \diam[z_2,x]) \ge
(C_1/2) (\diam A)^{-1} \diam[z_1,z_2],
$$
which concludes the proof of (\ref{eq-bt}). \qed

\medskip

{\em Proof of Corollary~\ref{cor-pcf}.}
In the p.c.f.\ case we have finitely many functions $f\in \Fk$, and each
of them
has an eventually periodic sequence of Taylor coefficients. Let $\ell$ be
so large that for any $f\ne h$ in $\Fk$ we have $|f\wedge h| < \ell$ and
for every $f\in \Fk$ the coefficients are periodic for $n\ge \ell$.
For any given $n\ge \ell$ there exists $L_n>0$ such that
we have
$$
(f\in \Fk,\ g\in \Bk,\ |f\wedge g| = n)\ \Rightarrow\ |g(\lam)| \ge
L_n |\lam|^n.
$$
Indeed, otherwise, passing to a subsequence we obtain $h\in \Bk$ with
$|f\wedge h|=n$ and $h(\lam)=0$, that is, $h\in \Fk$, contradicting the
choice of $n$ and $\ell$.

Let $p$ be a common period for all of the functions in $\Fk$. We claim that
$$
C_1 := \min\{L_n: \ \ell \le n \le \ell+p-1\}
$$
satisfies (\ref{eq-tech}). This is immediate from periodicity: if
$|f\wedge g| = n+kp$, then
$$
g(\lam) = g(\lam) - f(\lam) = \lam^{kp}[\wtil{g}(\lam) - f(\lam)] =
\lam^{kp}\wtil{g}(\lam)
$$
for some $\wtil{g}\in \Bk$ with $|f\wedge \wtil{g}|=n$, hence
we can take $L_{n+kp} = L_n$, and the claim follows. It remains to apply
Theorem~\ref{th-bt}. \qed


\begin{thebibliography}{99}


\bibitem{Ai} H. Aikawa, T. Lundh, and T. Mizutani,
Martin boundary of a fractal domain, {\em
Potential Anal.} {\bf 18} (2003), no.\ 4, 311--357.

\bibitem{bandt} C. Bandt, On the Mandelbrot set for pairs of linear maps,
{\em Nonlinearity} {\bf 15} (2002), 1127--1147.

\bibitem{bake} C. Bandt and K. Keller, Self-similar sets 2. A simple
approach to the topological structure of fractals, {\em Math.\ Nachr.}
{\bf 154} (1991), 27--39.

\bibitem{BR} C. Bandt and H. Rao,
Topology and separation of self-similar fractals in the plane,
{\em Nonlinearity} {\bf 20} (2007), 1463-1474.


\bibitem{BH} M. F. Barnsley and  A. N. Harrington, A Mandelbrot set
for pairs of linear maps, {\em Physica} {\bf 15D} (1985), 421--432.

\bibitem{BBBP} F. Beaucoup, P. Borwein, D. W. Boyd, and
C. Pinner, Multiple roots of $[-1,1]$ power series, {\em J.\ London Math.\
Soc.(2)} {\bf 57} (1998), 135--147.


\bibitem{bou2} T. Bousch, Connexit\'{e} locale et par chemins h\"{o}lderiens
pour les syst\`{e}mes it\'{e}r\'{e}s de fonctions, Preprint, 1993,
\verb+http://topo.math.u-psud.fr/~bousch+

\bibitem{BS} H. Bruin and D. Schleicher, Symbolic dynamics of quadratic polynomials, Preprint,
Institut Mittag-Leffler, 2002.

\bibitem{BS2} H. Bruin and D. Schleicher, Admissibility of kneading sequences and structure of Hubbard trees for quadratic polynomials,
Preprint arXiv:0801.4662v1, 2008.

\bibitem{CJY} L. Carleson, P. Jones, J.-Ch. Yoccoz, Julia and John.
{\em Bol. Soc. Bras. Mat.}
{\bf 25} (1994), 1--30.

\bibitem{Dev} R. L. Devaney, M. M. Rocha, S Siegmund, Rational maps with generalized Sierpi\'nski gasket Julia sets, {\em Topology and its Applications}
{\bf 154} (2007), 11--27.

\bibitem{DH} A. Douady and J. Hubbard, {\em Etudes Dynamiques des Polyn\^omes Complexes}, Publications Math\'ematique d'Orsay, 84-02, 1984, and 85-04, 1985.

\bibitem{GNV} M. Ghamsari, R. N\"akki and J. V\"ais\"al\"a, John disks and extension of maps, 
{\em Monatsh.\ Math.} {\bf 117} (1994), no.\ 1-2, 63--94.

\bibitem{HP} P. Ha\"issinsky, K. Pilgrim, Coarse expanding conformal dynamics, arxiv math.DS/0612617,
{\em Ast\'erisque}, to appear.

\bibitem{hata}
M. Hata, On the structure of self-similar sets,
{\em Japan.\ J.\ Appl.\ Math.} {\bf 2} (1985), 381--414.

\bibitem{H} J. Heinonen, Lectures on analysis on metric spaces. {\em Universitext. Springer-Verlag, New York}
(2001)

\bibitem{hutch} J.\ E.\ Hutchinson,
{\em Fractals and self-similarity.}
Indiana Univ.\ Math.\ J.\ 30 (1981), 713--747.

\bibitem{Kaffl} A. Kaffl, On the structure of abstract Hubbard trees and the space of abstract kneading 
sequences of degree two, {\em Ergodic Th.\ \& Dynam.\ Sys.} {\bf 27} (2007), 1215--1238.


\bibitem{kam1}
A. Kameyama, Julia sets and self-similar sets,
{\em  Topology Appl.} {\bf 54} (1993), 241--51.

\bibitem{kam2} A. Kameyama, Julia sets of post-critically finite
rational maps and topologically self-similar sets, {\em Nonlinearity}
{\bf 13} (2000), 165--188.

\bibitem{keller} K. Keller,
Invariant factors, Julia equivalences and the (abstract) Mandelbrot set,
{\em Lecture Notes in Mathematics} {\bf 1732}, Springer-Verlag, Berlin, 2000.


\bibitem{MS} C. McMullen and D. Sullivan,
Quasiconformal homeomorphisms and dynamics. III. The Teichm\"uller space of
a holomorphic dynamical system, {\em Adv.\ Math.} {\bf 135} (1998), 351--395.

\bibitem{Milnor} J. Milnor, {\em Dynamics in one complex variable}, Annals of Mathematics Studies, Princeton University Press, 2006.

\bibitem{Mih} N. Mihalache, {\em Julia and John revisited}, arXiv:0803.3889.


\bibitem{NV} R. N\"akki and J. V\"ais\"al\"a, John disks, {\em Exposition.\ Math.} {\bf 9} (1991), 3--43.

\bibitem{PR} F. Przytycki and S. Rohde,
Rigidity of holomorphic Collet-Eckmann repellers, {\em Ark.\ Math.} {\bf 37}
(1999), 357--371.

\bibitem{solo} B. Solomyak,
On the `Mandelbrot set' for pairs of linear maps: asymptotic self-similarity,
{\em  Nonlinearity}  {\bf 18}  (2005),  no.\ 5, 1927--1943.

\bibitem{solxu} B. Solomyak and H. Xu,
On the `Mandelbrot set' for a pair of linear maps and complex
Bernoulli convolutions, {\em Nonlinearity}
{\bf 16} (2003),
1733--1749.

\bibitem{TV} P. Tukia, J. V\"ais\"al\"a, Quasisymmetric embeddings of metric spaces.
{\em Ann. Acad. Sci. Fenn.} {\bf 5} (1980), 97--114.



\end{thebibliography}
\end{document}